\documentstyle[12pt,amssym]{article}
\newcommand{\A}{\mbox{\boldmath $A$}}
\newcommand{\C}{\mbox{\boldmath $C$}}
\newcommand{\F}{\mbox{\boldmath $F$}}
\newcommand{\I}{\mbox{\boldmath $I$}}
\newcommand{\g}{\mbox{\boldmath $g$}}
\newcommand{\B}{\mbox{\boldmath $B$}}
\newcommand{\D}{\mbox{\boldmath $D$}}
\newcommand{\BB}{\mbox{\boldmath $\cal B$}}

\def\CC{\mbox{$\Bbb C$}}

\def\case#1#2{{\textstyle{#1\over #2}}}

\newcommand{\tA}{\mbox{\hskip 3pt $\tilde{\mbox{\hskip -3pt $\A$}}$}}
\newcommand{\tB}{\mbox{\hskip 3pt $\tilde{\mbox{\hskip -3pt $\B$}}$}}
\newcommand{\tD}{\mbox{\hskip 3pt $\tilde{\mbox{\hskip -3pt $\D$}}$}}
\newcommand{\tBB}{\mbox{\hskip 3pt $\tilde{\mbox{\hskip -3pt $\BB$}}$}}

\begin{document}
\baselineskip=20pt plus 1pt minus1pt
\hfill ULB/229/CQ/98/8

\vspace{1.5cm}
\begin{center}

{\Large\bf\boldmath Covariant ($hh'$)-Deformed Bosonic and Fermionic Algebras
as Contraction Limits of $q$-Deformed Ones}\\

\bigskip

{\large\bf C.\ Quesne\footnote{\rm Physique Nucl\'eaire Th\'eorique et Physique
Math\'ematique, Universit\'e Libre de Bruxelles, B-1050 Brussels, Belgium}}

\end{center}

\vspace{10cm}
\noindent
Running head: Covariant ($hh'$)-Deformed Bosonic and Fermionic Algebras

\noindent
Mailing address: C. Quesne, Physique Nucl\'eaire Th\'eorique et Physique
Math\'ematique, Universit\'e Libre de Bruxelles, Campus de la Plaine CP229,
Boulevard du Triomphe, B-1050 Brussels, Belgium
\newpage
%
%
\begin{center}
{\bf Abstract}
\end{center}

\noindent
$\mbox{GL}_h(n) \times \mbox{GL}_{h'}(m)$-covariant ($hh'$)-bosonic (or
($hh'$)-fermionic) algebras ${\cal A}_{hh'\pm}(n,m)$ are built in terms of the
corresponding $R_h$ and $R_{h'}$-matrices by contracting the
$\mbox{GL}_q(n) \times \mbox{GL}_{q^{\pm1}}(m)$-covariant $q$-bosonic (or
$q$-fermionic) algebras ${\cal A}^{(\alpha)}_{q\pm}(n,m)$, $\alpha = 1$, 2. When
using a basis of ${\cal A}^{(\alpha)}_{q\pm}(n,m)$ wherein the annihilation
operators are contragredient to the creation ones, this contraction
procedure can
be carried out for any $n$, $m$ values. When employing instead a basis
wherein the
annihilation operators, as the creation ones, are irreducible tensor
operators with
respect to the dual quantum algebra $\mbox{U}_q(\mbox{gl}(n)) \otimes
\mbox{U}_{q^{\pm1}}(\mbox{gl}(m))$, a contraction limit only exists for
$n$, $m \in
\{1, 2, 4, 6,~\ldots\}$. For $n=2$, $m=1$, and $n=m=2$, the resulting
relations can
be expressed in terms of coupled (anti)commutators (as in the classical
case), by
using $\mbox{U}_h(\mbox{sl}(2))$ (instead of sl(2)) Clebsch-Gordan coefficients.
Some   U$_h$(sl(2)) rank-1/2 irreducible tensor operators, recently
constructed by
Aizawa, are shown to provide a realization of ${\cal A}_{h\pm}(2,1)$.
\newpage
%
%
\section{INTRODUCTION}
\label{sec:intro}
It is well known that the Lie group GL(2)  admits, up to isomorphism, only two
quantum group deformations with central determinant (Kupershmidt, 1992): the
standard deformation GL$_q$(2) (Drinfeld, 1987), and the so-called Jordanian
deformation GL$_h$(2) (Demidov {\em et al.}, 1990; Zakrzewski, 1991). On the
quantum algebra level, the Jordanian deformation U$_h$(sl(2)) of the classical
enveloping algebra U(sl(2)) was first considered by Ohn (1992), and its
universal ${\cal R}_h$-matrix was independently derived by Ballesteros and
Herranz (1996), and by Shariati {\em et al.}\ (1996). The fundamental
representation of U$_h$(sl(2)), which remains undeformed, was obtained by Ohn
(1992), while the other finite-dimensional highest-weight representations were
first studied by Dobrev (1996). Two-parametric Jordanian deformations
GL$_{h,\alpha}$(2), and U$_{h,\alpha}$(gl(2)) were also introduced by
Aghamohammadi (1993), Aneva {\em et al.}\ (1997), and Parashar (1998).\par
%
%
Two useful tools have been devised for studying the Jordanian deformations.
One of
them is a contraction procedure that allows one to construct the latter from
standard deformations (Aghamohammadi {\em et al.}, 1995): a similarity
transformation of the defining $R_q$ and $T_q$-matrices of GL$_q$(2) is
performed using a matrix singular itself in the $q \to 1$ limit, but in
such a way
that the transformed matrices are nonsingular, and yield the defining $R_h$ and
$T_h$-matrices of GL$_h$(2).\par
%
%
Such a contraction technique can be generalized to higher-dimensional quantum
groups. It was indeed shown by Alishahiha (1995) that there exist just two
independent singular maps from GL$_q$(3) to new quantum groups, one trivial and
one nontrivial, and that the latter can be extended to GL$_q$($N$) and
SP$_q$($2N$)
for arbitrary $N$. This gives rise to GL$_h$($N$) and SP$_h$($2N$),
respectively,
which are defined by their corresponding $R_h$-matrix.\par
%
%
The other tool consists in a class of nonlinear invertible maps between the
generators of U$_h$(sl(2)) and U(sl(2)) (Abdesselam {\em et al.}, 1998b).
Although
there exists an equivalence relation between these maps, they may arise
naturally
in different contexts, and may be particularly useful for different
purposes. One of
them (Abdesselam {\em et al.}, 1996) yields an explicit and simple
method for constructing the finite-dimensional irreducible
representations~(irreps)
of U$_h$(sl(2)). Furthermore, it provides the decomposition rule for the tensor
product of two such irreps (Aizawa, 1997), an explicit formula for U$_h$(sl(2))
Clebsch-Gordan coefficients (CGC) (Van der Jeugt, 1998), as well as bosonic and
fermionic realizations of irreducible tensor operators~(ITO) for U$_h$(sl(2)), and
an extension of Wigner-Eckart theorem to the latter (Aizawa, 1998). Another map
(Abdesselam {\em et al.}, 1998a) provides an operational generalization of the
contraction method of Aghamohammadi {\em et al.}\ (1995), and leads to the
construction of $R_h^{j_1;j_2}$ and $T_h^j$-matrices of arbitrary ($j_1 \otimes
j_2$) and $j$ irreps of U$_h$(sl(2)), respectively, as well as their
two-parametric
and/or coloured extensions (Chakrabarti and Quesne, 1998). Such a technique has
also been generalized to higher-dimensional quantum
algebras (Abdesselam {\em et al.}, 1997; Abdesselam {\em et al.}, 1998a).\par
%
%
In the present paper, we will apply the contraction procedure used by
Alishahiha (1995) to the $\mbox{GL}_q(n) \times \mbox{GL}_q(m)$-covariant
$q$-bosonic algebras ${\cal A}^{(\alpha)}_{q+}(n,m)$, $\alpha=1$, 2, and the
$\mbox{GL}_q(n) \times \mbox{GL}_{q^{-1}}(m)$-covariant $q$-fermionic ones
${\cal A}^{(\alpha)}_{q-}(n,m)$, which were constructed some years ago by the
present author (Quesne, 1993; Quesne, 1994), and recently rederived by
Fiore (1998) by another procedure. Such algebras generalize  Pusz-Woronowicz
GL$_q$($n$)-covariant $q$-bosonic or $q$-fermionic algebras (Pusz and
Woronowicz, 1989; Pusz, 1989), ${\cal A}^{(\alpha)}_{q\pm}(n)$, $\alpha=1$, 2,
to a tensor product of $m$ Fock spaces. They are generated by $nm$ pairs of
boson
or fermion-like creation and annihilation operators $\A^{\prime\dagger}_{is}$,
$\A'_{is}$ (or $\tA^{\prime\dagger}_{is}$), $i=1$, 2, $\ldots$,~$n$, $s=1$, 2,
$\ldots$,~$m$, with definite transformation properties under both GL$_q$($n$)
and GL$_{q^{\pm1}}$($m$), or U$_q$(gl($n$)) and U$_{q^{\pm1}}$(gl($m$)).\par
%
%
Our purpose will be twofold. Firstly, we will study under which conditions,
if any,
contracting these algebras by using two independent similarity
transformations for
GL$_q$($n$) and GL$_{q^{\pm1}}$($m$) may lead to $\mbox{GL}_h(n) \times
\mbox{GL}_{h'}(m)$-covariant ($hh'$)-bosonic or ($hh'$)-fermionic algebras
${\cal A}_{hh'\pm}(n,m)$. Secondly, in the $n=2$, $m=1$, and $n=m=2$ cases, we
will establish some relations with the works of Aizawa (1998) on ITO, and of Van
der Jeugt (1998) on CGC for U$_h$(sl(2)).\par
%
%
The algebras ${\cal A}_{hh'\pm}(n,m)$, whose generators $\A^+_{is}$,
$\A_{is}$ (or
$\tA_{is}$), $i=1$, 2, \ldots,~$n$, $s=1$, 2, \ldots,~$m$, have definite
transformation properties under both GL$_h$($n$) and GL$_{h'}$($m$), may be
useful in applications of Jordanian quantum groups in various fields, such as
quantum mechanics, condensed matter physics or quantum field theory. In such
applications, GL$_h$($n$) may represent the symmetry of the physical system,
while index $s$ may label different particles, crystal sites or space-time
points,
respectively. The deformed (anti)commutation relations satisfied by
$\A^+_{is}$, $\A_{is}$ (or $\tA_{is})$ may then either reflect some exotic
statistics or be interpreted as those of composite operators creating and
annihilating some quasi-particles or dressed states of bosons (or fermions).\par
%
%
This paper is organized as follows. Alishahiha's contraction procedure for
GL$_h$($N$) is reviewed in Sec.~\ref{sec:contraction}, and various forms of
$\mbox{GL}_q(n) \times \mbox{GL}_{q^{\pm1}}(m)$-covariant $q$-bosonic (or
$q$-fermionic) algebras are presented in Sec.~\ref{sec:covariant-q}. In
Sec.~\ref{sec:covariant-hh'}, the technique of Sec.~\ref{sec:contraction}
is applied
to such algebras to obtain $\mbox{GL}_h(n) \times \mbox{GL}_{h'}(m)$-covariant
($hh'$)-bosonic (or ($hh'$)-fermionic) algebras. The special cases where
$n=2$, and
$m=1$ or 2 are dealt with in Sec.~\ref{sec:special}.
Section~\ref{sec:conclusion}
contains the conclusion.\par
%
%
\section{\boldmath CONTRACTION OF GL$_q$($N$)}
\label{sec:contraction}
The quantum group GL$_q$($N$) is defined (Majid, 1990) as the associative
algebra
over \CC\ generated by $I$ and the noncommutative elements $T'_{ij}$ of an
$N\times N$ matrix $T'$ subject to the relations
\begin{equation}
  R'_q T'_1 T'_2 = T'_2 T'_1 R'_q, \qquad T'_1 = T' \otimes I, \qquad T'_2
= I \otimes
  T',    \label{eq:GL_q-alg}
\end{equation}
where
\begin{equation}
  R'_q = q \sum_i e_{ii} \otimes e_{ii} + \sum_{i\ne j} e_{ii} \otimes e_{jj}
  + \left(q - q^{-1}\right) \sum_{i<j} e_{ij} \otimes e_{ji},  \label{eq:R'_q}
\end{equation}
with $i$, $j$ running over 1, 2, \ldots,~$N$, and $e_{ij}$ denoting the $N
\times
N$~matrix with entry~1 in row~$i$ and column~$j$, and zeros everywhere else.
It is equipped with a coproduct, a counit, and an antipode defined by
\begin{equation}
  \Delta(T') = T'_1\, \dot{\otimes}\, T'_2, \qquad \epsilon(T') = I, \qquad
S(T') =
  T^{\prime -1},   \label{eq:GL_q-coalg}
\end{equation}
respectively, where $\dot{\otimes}$ denotes tensor product together with matrix
multiplication. An equivalent form of the $RTT$-relations (\ref{eq:GL_q-alg}) is
obtained by replacing $R'_q$ by $\tau R_q^{\prime -1} \tau$, where $\tau$ is the
twist map, i.e., $\tau (a \otimes b) = b \otimes a$. Note that throughout this
paper, $q$-deformed objects will be denoted by primed quantities, whereas
unprimed ones will represent $h$-deformed objects.\par
%
%
Let us consider the similarity transformation (Aghamohammadi {\em et al.}, 1995;
Alishahiha, 1995)
\begin{equation}
  R''_q = \left(g^{-1} \otimes g^{-1}\right) R'_q (g \otimes g), \qquad T''
= g^{-1} T'
  g,
\end{equation}
where $g$ is the $N \times N$ matrix defined by
\begin{equation}
  g = \sum_i e_{ii} + \eta e_{1N}, \qquad \eta = h/(q-1).  \label{eq:g}
\end{equation}
Eqs.\ (\ref{eq:GL_q-alg}) and (\ref{eq:GL_q-coalg}) simply become
\begin{equation}
  R''_q T''_1 T''_2 = T''_2 T''_1 R''_q, \qquad \Delta(T'') = T''_1\,
\dot{\otimes}\,
  T''_2,  \qquad \epsilon(T'') = I, \qquad S(T'') = T^{\prime\prime -1}.
  \label{eq:GL_q-sim}
\end{equation}
\par
%
%
When $q$ goes to one, although  parameter $\eta$ in (\ref{eq:g}) becomes
singular,
the relations in (\ref{eq:GL_q-sim}) have a definite limit
\begin{equation}
  R_h T_1 T_2 = T_2 T_1 R_h, \qquad \Delta(T) = T_1\, \dot{\otimes}\, T_2,
\qquad
  \epsilon(T) = I, \qquad S(T) = T^{-1},  \label{eq:GL_h}
\end{equation}
where $T \equiv \lim_{q\to1} T''$, and
\begin{eqnarray}
  R_h & \equiv & \lim_{q\to1} R''_q \nonumber \\
  & = & \sum_{ij} e_{ii} \otimes e_{jj} + h \biggl[e_{11} \otimes e_{1N} -
e_{1N}
       \otimes e_{11} + e_{1N} \otimes e_{NN} - e_{NN} \otimes e_{1N}
\nonumber \\
  & & \mbox{} + 2 \sum_{i=2}^{N-1} (e_{1i} \otimes e_{iN} - e_{iN} \otimes
e_{1i})
       \biggr] + h^2 e_{1N} \otimes e_{1N}.  \label{eq:R_h}
\end{eqnarray}
The resulting $R_h$-matrix is triangular, i.e., it is quasitriangular and
$R_h = \tau
R_h^{-1} \tau$, showing that the two equivalent forms of $RTT$-relations for
GL$_q$($N$) have actually the same contraction limit. Together with $I$, the
elements $T_{ij}$ of the $N \times N$ matrix $T$ generate the Jordanian quantum
group GL$_h$($N$).\par
%
%
\section{\boldmath COVARIANT $q$-BOSONIC AND $q$-FERMIONIC ALGEBRAS}
\label{sec:covariant-q}
\setcounter{equation}{0}
Let us consider two different copies of the quantum group GL$_q$($N$) considered
in Sec.~\ref{sec:contraction}, corresponding to possibly different dimensions
$n$, $m$, and parameters $q$, $q^{\sigma}$, respectively. Let us denote
quantities referring  to GL$_q$($n$) by ordinary (primed) letters ($R'_q$,
$T'$,~\ldots), and quantities referring to GL$_{q^{\sigma}}$($m$) by script
(primed)
letters (${\cal R}'_{q^{\sigma}}$, $\cal{T}'$,~\ldots). The elements
$T'_{ij}$, $i$,
$j=1$, 2,~$\ldots n$, of GL$_q$($n$) are assumed to commute with the elements
${\cal T}'_{st}$, $s$, $t=1$, 2,~$\ldots m$, of GL$_{q^{\sigma}}$($m$).
Note that for
simplicity's sake, we have skipped the parameters $q$ and $q^{\sigma}$, which
should be appended to $T'$ and $\cal{T}'$, respectively. With GL$_q$($n$) and
GL$_{q^{\sigma}}$($m$), we can associate the dual (commuting) quantum algebras
U$_q$(gl($n$)) and U$_{q^{\sigma}}$(gl($m$)).\par
%
%
Some years ago, it was shown (Quesne, 1993) that $q$-bosonic creation  and
annihilation operators $\A^{\prime +}_{is}$, $\tA'_{is}$, $i=1$, 2,
$\ldots$,~$n$, $s=1$, 2, $\ldots$,~$m$, that are double ITO of rank $[1
\dot{0}]_n
\times [1 \dot{0}]_m$ and $[\dot{0} -1]_n \times [\dot{0} -1]_m$  with
respect to
the quantum algebra $\mbox{U}_q(\mbox{gl}(n)) \times \mbox{U}_q(\mbox{gl}(m))$,
respectively, can be constructed in terms of standard $q$-bosonic creation,
annihilation, and number operators $a^{\prime +}_{is}$, $a'_{is}$,
$N'_{is}$, $i=1$,
2, $\ldots$,~$n$, $s=1$, 2, $\ldots$,~$m$ (Biedenharn, 1989; Macfarlane, 1989;
Sun and Fu, 1989), acting in a tensor product Fock space $\F = \prod_{i=1}^n
\prod_{s=1}^m \otimes F_{is}$. Here $[1 \dot{0}]_n$ and $[\dot{0} -1]_n$ denote
$n$-row Young diagrams, the dot over 0 meaning that this numeral is repeated as
often as necessary. It is straightforward to extend such a
construction to covariant $q$-fermionic operators, provided one replaces
$\mbox{U}_q(\mbox{gl}(m))$ by $\mbox{U}_{q^{-1}}(\mbox{gl}(m))$, and standard
$q$-bosonic operators by standard $q$-fermionic ones (Chaichian and Kulish,
1990; Hayashi, 1990).\par
%
%
The annihilation operators $\A'_{is}$, contragredient to
$\A^{\prime+}_{is}$, can
also be considered, and are related to the covariant ones $\tA'_{is}$
through the
equation
\begin{equation}
  \tA'_{is} = (-1)^{i+s} q^{[n-2i+1 + \sigma(m-2s+1)]/2} \A'_{i's'},
  \label{eq:A'-tilde}
\end{equation}
where $i' \equiv n+1-i$, $s' \equiv m+1-s$, and $\sigma = +1$ (resp.\ $-1$) for
$q$-bosons (resp.\ $q$-fermions). In matrix form, Eq.~(\ref{eq:A'-tilde}) can be
rewritten as
\begin{equation}
  \tA' = \A' \C', \qquad \C' = C'_q {\cal C}'_{q^{\sigma}},
\end{equation}
where
\begin{equation}
  C'_q = \sum_i (-1)^{n-i} q^{-(n-2i+1)/2} e_{ii'}, \qquad {\cal
C}'_{q^{\sigma}} =
  \sum_s (-1)^{m-s} q^{-\sigma(m-2s+1)/2} e_{ss'}.
\end{equation}
\par
%
%
As it happens in the $m=1$ case for the GL$_q$($n$)-covariant $q$-bosonic or
$q$-fermionic operators (Pusz and Woronowicz, 1989; Pusz, 1989), there actually
exist two independent ways of constructing $\A^{\prime+}_{is}$ and
$\tA'_{is}$ (or $\A'_{is}$) in terms of $a^{\prime +}_{is}$, $a'_{is}$,
$N'_{is}$. According to the choice made, the operators $\A^{\prime+}_{is}$ and
$\tA'_{is}$, or $\A^{\prime+}_{is}$ and $\A'_{is}$, generate with $\I = I
{\cal I}$ one of two different $\mbox{U}_q(\mbox{gl}(n)) \times
\mbox{U}_{q^{\sigma}}(\mbox{gl}(m))$-module, or $\mbox{GL}_q(n) \times
\mbox{GL}_{q^{\sigma}}(m)$-comodule algebras, which will be denoted by ${\cal
A}^{(1)}_{q\sigma}(n,m)$ and ${\cal A}^{(2)}_{q\sigma}(n,m)$. The defining
relations
of such algebras can be written in two compact forms, enhancing the
transformation properties of the operators under the quantum group
$\mbox{GL}_q(n) \times \mbox{GL}_{q^{\sigma}}(m)$ or the corresponding quantum
algebra $\mbox{U}_q(\mbox{gl}(n)) \times \mbox{U}_{q^{\sigma}}(\mbox{gl}(m))$,
respectively, as well as in componentwise form using $q$-(anti)commutators.\par
%
%
In the first compact form, the defining relations of ${\cal
A}^{(1)}_{q\sigma}(n,m)$ in the $\left\{\A^{\prime+}_{is}, \A'_{is} \right\}$
basis read (Quesne, 1994; Fiore, 1998)
\begin{equation}
  R'_q \A^{\prime +}_1 \A^{\prime +}_2 = \sigma \A^{\prime +}_2 \A^{\prime +}_1
  {\cal R}'_{q^{\sigma}},  \label{eq:compact1}
\end{equation}
\begin{equation}
  R'_q \A'_2 \A'_1 = \sigma \A'_1 \A'_2 {\cal R}'_{q^{\sigma}},
  \label{eq:compact2}
\end{equation}
\begin{equation}
  \A'_2 \A^{\prime +}_1 = \I_{21} + \sigma R^{\prime t_1}_q
  {\cal R}^{\prime t_1}_{q^{\sigma}} \A^{\prime +}_1 \A'_2,  \label{eq:compact3}
\end{equation}
while those of ${\cal A}^{(2)}_{q\sigma}(n,m)$ are given by
Eqs.~(\ref{eq:compact1}), (\ref{eq:compact2}), and
\begin{equation}
  \A'_1 \A^{\prime +}_2 = \I_{12} + \sigma R^{\prime t_2}_{q^{-1}}
  {\cal R}^{\prime t_2}_{q^{-\sigma}} \A^{\prime +}_2 \A'_1.
  \label{eq:compact3bis}
\end{equation}
Here we use the defining $R'_q$-matrix of GL$_q$($n$), given in Eq.\
(\ref{eq:R'_q}),
and its counterpart ${\cal R}'_{q^{\sigma}}$ for GL$_{q^{\sigma}}$($m$), as
well as a
shorthand tensor notation similar to that of Eq.\ (\ref{eq:GL_q-alg}), with
$t_1$
(resp.\ $t_2$) denoting transposition in the first (resp.\ second) space of
the tensor
product.\par
%
%
When using instead the $\left\{\A^{\prime +}_{is}, \tA'_{is}\right\}$ basis of
${\cal A}^{(1)}_{q\sigma}(n,m)$ and ${\cal A}^{(2)}_{q\sigma}(n,m)$,
Eqs.\ (\ref{eq:compact2}), (\ref{eq:compact3}), and (\ref{eq:compact3bis})
become
(Quesne, 1994)
\begin{equation}
  R'_q \tA'_1 \tA'_2 = \sigma \tA'_2 \tA'_1 {\cal R}'_{q^{\sigma}},
  \label{eq:compact2-tilde}
\end{equation}
\begin{equation}
  \tA'_2 \A^{\prime +}_1 = \C'_{12} + \sigma \A^{\prime +}_1 \tA'_2
  \tilde{R}^{\prime -1}_q \tilde{{\cal R}}^{\prime -1}_{q^{\sigma}},
  \label{eq:compact3-tilde}
\end{equation}
and
\begin{equation}
  \tA'_1 \A^{\prime +}_2 = \C'_{21} + \sigma \A^{\prime +}_2 \tA'_1
  \tilde{R}'_q \tilde{{\cal R}}'_{q^{\sigma}},  \label{eq:compact3bis-tilde}
\end{equation}
where
\begin{equation}
  \tilde{R}'_q \equiv C^{\prime -1}_{q,1} \left(R^{\prime
-1}_q\right)^{t_1} C'_{q,1}
  = C^{\prime -1}_{q,2} \left(R^{\prime t_2}_q\right)^{-1} C'_{q,2},
\end{equation}
and similarly for $\tilde{{\cal R}}'_{q^{\sigma}}$. Note that one can go
from ${\cal
A}^{(1)}_{q\sigma}(n,m)$ to ${\cal A}^{(2)}_{q\sigma}(n,m)$ by making the
substitutions $R'_q \to \tau R^{\prime -1}_q \tau$, ${\cal R}'_{q^{\sigma}}
\to \tau
{\cal R}^{\prime -1}_{q^{\sigma}} \tau$.\par
%
%
In either form (\ref{eq:compact1})--(\ref{eq:compact3}) (resp.\
(\ref{eq:compact1}), (\ref{eq:compact2}), (\ref{eq:compact3bis})) or
(\ref{eq:compact1}), (\ref{eq:compact2-tilde}), (\ref{eq:compact3-tilde})
(resp.\
(\ref{eq:compact1}), (\ref{eq:compact2-tilde}),
(\ref{eq:compact3bis-tilde})), it is
easy to see that ${\cal A}^{(1)}_{q\sigma}(n,m)$ (resp.\ ${\cal A}^{
(2)}_{q\sigma}(n,m)$) is a $\mbox{GL}_q(n) \times
\mbox{GL}_{q^{\sigma}}(m)$-comodule algebra. The transformation
\begin{equation}
  \varphi'\left(\A^{\prime +}\right) = \A^{\prime +} T' {\cal T}', \qquad
  \varphi'\left(\A'\right) = T^{\prime -1} {\cal T}^{\prime -1} \A',
\label{eq:phi'}
\end{equation}
or
\begin{equation}
  \varphi'\left(\A^{\prime +}\right) = \A^{\prime +} T' {\cal T}', \qquad
  \varphi'\left(\tA'\right) = \tA' \tilde{T}' \tilde{{\cal T}}',
\end{equation}
where $T'_{ij} \in \mbox{GL}_q(n)$, ${\cal T}'_{st} \in
\mbox{GL}_{q^{\sigma}}(m)$,
$\tilde{T}' = C^{\prime -1}_q \left(T^{\prime -1}\right)^t C'_q$, and
$\tilde{{\cal T}}'
= {\cal C}^{\prime -1}_{q^{\sigma}} \left({\cal T}^{\prime -1}\right)^t {\cal
C}'_{q^{\sigma}}$, indeed leaves the defining equations invariant, while being
consistent with the $\mbox{GL}_q(n) \times \mbox{GL}_{q^{\sigma}}(m)$ coalgebra
structure, as given in Eq.\ (\ref{eq:GL_q-coalg}), and its counterpart for
GL$_{q^{\sigma}}$($m$).\par
%
%
{}For $m=1$, one gets ${\cal R}'_{q^{\sigma}} = q^{\sigma}$, ${\cal
C}'_{q^{\sigma}} =
1$, $\tilde{{\cal R}}'_{q^{\sigma}} = q^{-\sigma}$, so that the defining
relations of
${\cal A}^{(1)}_{q\sigma}(n,1)$ and ${\cal A}^{(2)}_{q\sigma}(n,1)$
coincide with
those of the two independent Pusz-Woronowicz algebras (Pusz and Woronowicz,
1989; Pusz, 1989).\par
%
%
The second compact form uses coupled $q$-(anti)commutators, defined by (Quesne,
1993)
\begin{eqnarray}
  & &\left[T^{[\lambda_1]_n [\lambda'_1]_m}, U^{[\lambda_2]_n
          [\lambda'_2]_m}\right\}^{[\Lambda]_n [\Lambda']_m}_{(M)_n (M')_m
          q^{\alpha}} = \left[T^{[\lambda_1]_n [\lambda'_1]_m} \times
U^{[\lambda_2]_n
          [\lambda'_2]_m}\right]^{[\Lambda]_n [\Lambda']_m}_{(M)_n (M')_m}
          \nonumber \\
  & & \qquad \mbox{} - \sigma (-1)^{\epsilon} q^{\alpha} \left[U^{[\lambda_2]_n
          [\lambda'_2]_m} \times T^{[\lambda_1]_n [\lambda'_1]_m}\right]
          ^{[\Lambda]_n [\Lambda']_m}_{(M)_n (M')_m}.  \label{eq:coupledcom}
\end{eqnarray}
Here the left-hand side is a coupled $q$-commutator (resp.\ $q$-anticommutator)
for $\sigma = +1$ (resp.\ $-1$), $T^{[\lambda_1]_n [\lambda'_1]_m}$ and
$U^{[\lambda_2]_n [\lambda'_2]_m}$ denote two double ITO of rank $[\lambda_1]_n
\times [\lambda'_1]_m$ and $[\lambda_2]_n \times [\lambda'_2]_m$ with respect
to $\mbox{U}_q(\mbox{gl}(n)) \times \mbox{U}_{q^{\sigma}}(\mbox{gl}(m))$,
respectively, their tensor product of rank $[\Lambda]_n \times [\Lambda']_m$ is
defined by
\begin{eqnarray}
  &&\left[T^{[\lambda_1]_n [\lambda'_1]_m} \times U^{[\lambda_2]_n
          [\lambda'_2]_m}\right]^{[\Lambda]_n [\Lambda']_m}_{(M)_n (M')_m} =
          \sum_{(\mu_1)_n (\mu'_1)_m (\mu_2)_n (\mu'_2)_m}
          \left\langle [\lambda_1]_n (\mu_1)_n, [\lambda_2]_n (\mu_2)_n |
          [\Lambda]_n (M)_n\right\rangle_q \nonumber \\
  && \mbox{} \times \left\langle [\lambda'_1]_m (\mu'_1)_m,
          [\lambda'_2]_m (\mu'_2)_m |  [\Lambda']_m (M')_m
          \right\rangle_{q^{\sigma}} T^{[\lambda_1]_n [\lambda'_1]_m}_{(\mu_1)_n
          (\mu'_1)_m} U^{[\lambda_2]_n [\lambda'_2]_m}_{(\mu_2)_n (\mu'_2)_m},
          \label{eq:tensorprod}
\end{eqnarray}
and the phase factor $\epsilon$ is given by
\begin{equation}
  \epsilon = \phi([\lambda_1]_n) + \phi([\lambda_2]_n) - \phi([\Lambda]_n) +
  \phi([\lambda'_1]_m)  + \phi([\lambda'_2]_m) - \phi([\Lambda']_m),
  \label{eq:epsilon}
\end{equation}
\begin{equation}
  \phi([\lambda_1]_n) = \case{1}{2} \sum_{i=1}^n (n+1-2i) \lambda_{1i}, \qquad
  \phi([\lambda'_1]_m) = \case{1}{2} \sum_{s=1}^m (m+1-2s) \lambda'_{1s}.
  \label{eq:phi-lambda}
\end{equation}
In Eq.\ (\ref{eq:tensorprod}), $\langle\, , | \,\rangle_q$ and $\langle\, , |
\,\rangle_{q^{\sigma}}$ denote U$_q$(gl($n$)) and U$_{q^{\sigma}}$(gl($m$))
CGC (Biedenharn, 1990), respectively, and we have assumed that the couplings are
multiplicity free (which is the case for the generators of ${\cal
A}^{(1)}_{q\sigma}(n,m)$ and ${\cal A}^{(2)}_{q\sigma}(n,m)$).\par
%
%
Such a compact form only exists for the $\{\A^{\prime+}_{is}, \tA'_{is}\}$
basis,
since $\A^{\prime+}_{is}$ and $\tA'_{is}$ (but not $\A'_{is}$) have a
definite rank
with respect to $\mbox{U}_q(\mbox{gl}(n)) \times
\mbox{U}_{q^{\sigma}}(\mbox{gl}(m))$, namely $[1 \dot{0}]_n \times [1\dot{0}]_m$
and $[\dot{0} -1]_n \times [\dot{0} -1]_m$, respectively. For
${\cal A}^{(1)}_{q\sigma}(n,m)$, one finds (Quesne, 1993)
\begin{equation}
  \left[\A^{\prime +}, \A^{\prime +}\right]^{[2\dot{0}]_n [1^2\dot{0}]_m}
  =  \left[\A^{\prime +}, \A^{\prime +}\right]^{[1^2\dot{0}]_n
  [2\dot{0}]_m} = 0,  \label{eq:coupledcom1-1}
\end{equation}
\begin{equation}
  \left[\tA', \tA'\right]^{[\dot{0}-2]_n [\dot{0}(-1)^2]_m} = \left[\tA',
  \tA'\right]^{[\dot{0}(-1)^2]_n [\dot{0}-2]_m} = 0,  \label{eq:coupledcom1-2}
\end{equation}
in the $q$-bosonic case ($\sigma = +1$), or
\begin{equation}
  \left\{\A^{\prime +}, \A^{\prime +}\right\}^{[2\dot{0}]_n [2\dot{0}]_m}
  = \left\{\A^{\prime +}, \A^{\prime +}\right\}^{[1^2\dot{0}]_n
[1^2\dot{0}]_m} = 0,
\end{equation}
\begin{equation}
  \left\{\tA', \tA'\right\}^{[\dot{0}-2]_n [\dot{0}-2]_m} = \left\{\tA',
  \tA'\right\}^{[\dot{0}(-1)^2]_n [\dot{0}(-1)^2]_m} = 0,
\end{equation}
in the $q$-fermionic one ($\sigma = -1$), and
\begin{equation}
  \left[\tA', \A^{\prime +}\right\}^{[1\dot{0}-1]_n [1\dot{0}-1]_m} =
\left[\tA',
  \A^{\prime+}\right\}^{[1\dot{0}-1]_n [\dot{0}]_m}_{q^{\sigma m}} = \left[\tA',
  \A^{\prime+}\right\}^{[\dot{0}]_n [1\dot{0}-1]_m}_{q^n} = 0,
  \label{eq:coupledcom3-1}
\end{equation}
\begin{equation}
  \left[\tA', \A^{\prime +}\right\}^{[\dot{0}]_n [\dot{0}]_m}
       _{q^{n+\sigma m}} = \sqrt{[n]_q [m]_q} \, \I, \label{eq:coupledcom3-2}
\end{equation}
in both cases ($\sigma = \pm 1$). For simplicity's sake, we have not written the
$\mbox{U}_q(\mbox{gl}(n)) \times \mbox{U}_{q^{\sigma}}(\mbox{gl}(m))$ irrep
row labels $(M_1)_n (M_2)_m$. As usual, $q$-numbers are defined by $[x]_q \equiv
\left(q^x - q^{-x}\right)/\left(q - q^{-1}\right)$. For ${\cal
A}^{(2)}_{q\sigma}(n,m)$, Eqs.\
(\ref{eq:coupledcom1-1})--(\ref{eq:coupledcom3-2}) remain valid but for the
substitution $q\to q^{-1}$ in the lower subscripts in Eqs.\
(\ref{eq:coupledcom3-1})
and (\ref{eq:coupledcom3-2}).\par
%
%
By using the explicit form of the $R'_q$ and ${\cal R}'_{q^{\sigma}}$ matrix
elements given in Eq.~(\ref{eq:R'_q}), or the explicit values of the
U$_q$(gl($n$))
and U$_{q^{\sigma}}$(gl($m$)) CGC (Biedenharn, 1990) together with
Eq.~(\ref{eq:A'-tilde}), Eqs.\ (\ref{eq:compact1})--(\ref{eq:compact3}),
or (\ref{eq:coupledcom1-1})--(\ref{eq:coupledcom3-2}), for ${\cal
A}^{(1)}_{q\sigma}(n,m)$ can be rewritten in componentwise form. The
results read (Quesne, 1993)
\begin{equation}
  \left\{\A^{\prime +}_{is}, \A^{\prime +}_{is}\right\} = 0,  \label{eq:com1}
\end{equation}
in the $q$-fermionic case ($\sigma = -1$), and
\begin{equation}
  \left[\A^{\prime +}_{is}, \A^{\prime +}_{it}\right\}_{q^{-1}}= 0, \qquad s<t,
  \label{eq:com2}
\end{equation}
\begin{equation}
  \left[\A^{\prime +}_{is}, \A^{\prime +}_{js}\right\}_{q^{-\sigma}} = 0,
\qquad i<j,
\end{equation}
\begin{equation}
  \left[\A^{\prime +}_{is}, \A^{\prime +}_{jt}\right\} = 0, \qquad i>j, s<t,
\end{equation}
\begin{equation}
  \left[\A^{\prime +}_{is}, \A^{\prime +}_{jt}\right\} = - \left(q -
q^{-1}\right)
  \A^{\prime +}_{js} \A^{\prime +}_{it}, \qquad i<j, s<t,  \label{eq:com3}
\end{equation}
\begin{equation}
  \left[\A'_{is}, \A^{\prime +}_{jt}\right\} = 0, \qquad i\ne j, s\ne t,
  \label{eq:com4}
\end{equation}
\begin{equation}
  \left[\A'_{is}, \A^{\prime +}_{js}\right\}_{q^{\sigma}} = \left(q -
q^{-1}\right)
  \sum_{t=1}^{s-1} \A^{\prime +}_{jt} \A'_{it},  \qquad i\ne j, \label{eq:com5}
\end{equation}
\begin{equation}
  \left[\A'_{is}, \A^{\prime +}_{it}\right\}_q = \sigma \left(q - q^{-1}\right)
  \sum_{j=1}^{i-1} \A^{\prime +}_{jt} \A'_{js}, \qquad  s\ne t,
\end{equation}
\begin{eqnarray}
  \left[\A'_{is}, \A^{\prime +}_{is}\right\}_{q^{1+\sigma}} & = & \I +
          \left(q^{2\sigma} - 1\right) \sum_{j=1}^{i-1} \A^{\prime +}_{js}
\A'_{js} +
          \left(q^2 - 1\right) \sum_{t=1}^{s-1} \A^{\prime +}_{it} \A'_{it}
\nonumber \\
  & & + \left(q - q^{-1}\right)^2\, \sum_{j=1}^{i-1} \sum_{t=1}^{s-1}
          \A^{\prime +}_{jt} \A'_{jt}, \label{eq:com6}
\end{eqnarray}
in both $q$-bosonic and $q$-fermionic cases ($\sigma = \pm 1$), together
with the
Hermitian conjugates of Eqs.\ (\ref{eq:com1})--(\ref{eq:com3}) (for real
$q$). Here,
for $q$-bosons (resp.\ $q$-fermions), $[\,,\,\}$ denotes a commutator (resp.\
anticommutator), and $[\,,\,\}_{q^{\alpha}}$ a $q$-commutator (resp.\
$q$-anticommutator), i.e., $[A, B\}_{q^{\alpha}} \equiv AB - \sigma q^{\alpha}
BA$.\par
%
%
{}For ${\cal A}^{(2)}_{q\sigma}(n,m)$, Eqs.\ (\ref{eq:com1})--(\ref{eq:com4})
remain unchanged, whereas Eqs.\ (\ref{eq:com5})--(\ref{eq:com6}) are replaced by
\begin{equation}
  \left[\A'_{is}, \A^{\prime +}_{js}\right\}_{q^{-\sigma}} = - \left(q -
q^{-1}\right)
  \sum_{t=s+1}^{m} \A^{\prime +}_{jt} \A'_{it},  \qquad i\ne j,
\end{equation}
\begin{equation}
  \left[\A'_{is}, \A^{\prime +}_{it}\right\}_{q^{-1}} = - \sigma \left(q -
  q^{-1}\right) \sum_{j=i+1}^{n} \A^{\prime +}_{jt} \A'_{js}, \qquad  s\ne t,
\end{equation}
\begin{eqnarray}
  \left[\A'_{is}, \A^{\prime +}_{is}\right\}_{q^{-1-\sigma}} & = & \I +
          \left(q^{-2\sigma} - 1\right) \sum_{j=i+1}^{n} \A^{\prime +}_{js}
\A'_{js} +
          \left(q^{-2} - 1\right) \sum_{t=s+1}^{m} \A^{\prime +}_{it} \A'_{it}
          \nonumber \\
  & & \mbox{} + \left(q - q^{-1}\right)^2 \sum_{j=i+1}^{n} \sum_{t=s+1}^{m}
          \A^{\prime +}_{jt} \A'_{jt}.  \label{eq:com6-bis}
\end{eqnarray}
\par
%
%
Note again that for $m=1$, Eqs.~(\ref{eq:com1})--(\ref{eq:com6-bis}) give
back the
Pusz-Woronowicz results (Pusz and Woronowicz, 1989; Pusz, 1989).\par
%
%
\section{\boldmath COVARIANT ($hh'$)-BOSONIC AND ($hh'$)-FERMIONIC ALGEBRAS}
\label{sec:covariant-hh'}
\setcounter{equation}{0}
Let us apply the contraction procedure of Sec.\ \ref{sec:contraction} to the
$\mbox{GL}_q(n) \times \mbox{GL}_{q^{\sigma}}(m)$-covariant $q$-bosonic (or
$q$-fermionic) algebras ${\cal A}^{(1)}_{q\sigma}(n,m)$ and ${\cal
A}^{(2)}_{q\sigma}(n,m)$. We shall successively consider the cases where
they are defined in the $\left\{\A^{\prime +}_{is}, \A'_{is}\right\}$
basis, or in the
$\left\{\A^{\prime +}_{is}, \tA'_{is}\right\}$ one.\par
%
%
Since we now have two commuting copies of GL$_q$($N$), we have to consider two
transformation matrices of type (\ref{eq:g}), $g = \sum_i e_{ii} + \eta
e_{1n}$, and
$\mbox{\textsf g} = \sum_s e_{ss} + \eta' e_{1m}$. They act on GL$_q$($n$) and
GL$_{q^{\sigma}}$($m$), respectively, and depend upon two parameters $\eta
\equiv h/(q-1)$, and $\eta' \equiv h'/\left(q^{\sigma}-1\right)$, which we may
assume independent.\par
%
%
Let us first consider Eqs.\ (\ref{eq:compact1})--(\ref{eq:compact3}), defining
${\cal A}^{(1)}_{q\sigma}(n,m)$ in the $\left\{\A^{\prime +}_{is},
\A'_{is}\right\}$ basis, and introduce transformed $q$-bosonic (or
$q$-fermionic) operators $\A^{\prime\prime +} = \A^{\prime +} \g$, $\A'' =
\g^{-1}
\A'$, where $\g = g\, \mbox{\textsf g}$, i.e., $\g_{is,jt} = g_{ij}\,
\mbox{\textsf g}_{st}$. By using the property $R^{\prime t}_q = \tau R'_q \tau$,
satisfied by (\ref{eq:R'_q}), and a similar one for ${\cal
R}'_{q^{\sigma}}$, it is
straightforward to show that Eqs.\ (\ref{eq:compact1})--(\ref{eq:compact3})
become
\begin{equation}
  \A^{\prime\prime +}_1 \A^{\prime\prime +}_2 = \sigma
  \A^{\prime\prime +}_2 \A^{\prime\prime +}_1 \left(\tau R''_{q^{-1}}
  \tau\right) {\cal R}''_{q^{\sigma}},  \label{eq:compact1-transf}
\end{equation}
\begin{equation}
  \A''_1 \A''_2 = \sigma R''_q \left(\tau {\cal R}''_{q^{-\sigma}}
  \tau\right) \A''_2 \A''_1  ,
\end{equation}
\begin{equation}
  \A''_2 \A^{\prime\prime +}_1 = \I_{21} + \sigma
  R^{\prime\prime t_1}_q {\cal R}^{\prime\prime t_1}_{q^{\sigma}}
  \A^{\prime\prime +}_1 \A''_2.  \label{eq:compact3-transf}
\end{equation}
\par
%
%
Defining now ($hh'$)-bosonic (or ($hh'$)-fermionic) operators by
\begin{equation}
  \A^+_{is} \equiv \lim_{q\to1} \A^{\prime\prime +}_{is}, \qquad \A_{is} \equiv
  \lim_{q\to1} \A''_{is},
\end{equation}
and taking the $q\to1$ limit of Eqs.\
(\ref{eq:compact1-transf})--(\ref{eq:compact3-transf}), we obtain that together
with $\I$, they generate an algebra ${\cal A}_{hh'\sigma}(n,m)$, whose defining
relations are
\begin{equation}
  \A^+_1 \A^+_2 = \sigma \A^+_2 \A^+_1 R_h {\cal R}_{h'}, \label{eq:compact1-h}
\end{equation}
\begin{equation}
  \A_1 \A_2 = \sigma R_h {\cal R}_{h'} \A_2 \A_1, \label{eq:compact2-h}
\end{equation}
\begin{equation}
  \A_2 \A^+_1 = \I_{21} + \sigma R^{t_1}_h {\cal R}^{t_1}_{h'} \A^+_1 \A_2.
  \label{eq:compact3-h}
\end{equation}
In deriving the latter, we explicitly used the fact that both $R_h$ and ${\cal
R}_{h'}$ are triangular. Similarly, transformation (\ref{eq:phi'}) goes into
\begin{equation}
  \varphi\left(\A^+\right) = \A^+ T {\cal T}, \qquad
  \varphi\left(\A\right) = T^{-1} {\cal T}^{-1} \A,  \label{eq:phi}
\end{equation}
where $T_{ij} \in \mbox{GL}_h(n)$, ${\cal T}_{st} \in \mbox{GL}_{h'}(m)$, and
$\varphi$ leaves Eqs.\ (\ref{eq:compact1-h})--(\ref{eq:compact3-h}) invariant,
while being consistent with the $\mbox{GL}_h(n) \times \mbox{GL}_{h'}(m)$
coalgebra structure, as given by Eq.\ (\ref{eq:GL_h}). Hence, ${\cal
A}_{hh'\sigma}(n,m)$ is a $\mbox{GL}_h(n) \times \mbox{GL}_{h'}(m)$-covariant
($hh'$)-bosonic (or ($hh'$)-fermionic) algebra.\par
%
%
It is easy to see that the same procedure applied to Eqs.\ (\ref{eq:compact1}),
(\ref{eq:compact2}), and (\ref{eq:compact3bis}), defining ${\cal A}^{
(2)}_{q\sigma}(n,m)$ in the $\left\{\A^{\prime +}_{is}, \A'_{is}\right\}$ basis,
leads to the same equations (\ref{eq:compact1-h})--(\ref{eq:compact3-h}) because
$R_h$ and ${\cal R}_{h'}$ are triangular. The algebra ${\cal
A}_{hh'\sigma}(n,m)$ is
therefore the contraction limit of both ${\cal A}^{(1)}_{q\sigma}(n,m)$ and
${\cal A}^{(2)}_{q\sigma}(n,m)$.\par
%
%
{}From Eqs.\ (\ref{eq:compact1-h}) and (\ref{eq:compact2-h}), it is clear that
contrary to what happens in the $q$-deformed case, $\A_{is}$ can never be
considered as the adjoint of $\A^+_{is}$. This comes from the lack of
*-structure on GL$_h$($N$).\par
%
%
Equations~(\ref{eq:compact1-h})--(\ref{eq:compact3-h}) agree with the general
form of $\cal H$-covariant deformed bosonic (or fermionic) algebras for
triangular Hopf algebras $\cal H$, which was derived by Fiore (1997). In the
present paper, we did establish that they can be obtained in a straightforward
way by Alishahiha's contraction technique (Alishahiha, 1995).\par
%
%
By using the explicit expression of $R_h$, given in Eq.\ (\ref{eq:R_h}),
and a similar
one for ${\cal R}_{h'}$, Eqs.\ (\ref{eq:compact1-h})--(\ref{eq:compact3-h})
can be
rewritten in componentwise form as follows:
\begin{eqnarray}
  \left[\A^+_{is}, \A^+_{jt}\right\} & = & (1 - \sigma P_{ij} P_{st})
         \bigl\{h \delta_{j,n} \left(1 - \delta_{\sigma,-1} \delta_{i,1}
\delta_{s,t}
         \right) d_i \A^+_{1s} \A^+_{it} \nonumber \\
  && \mbox{} + h' \delta_{t,m} \left(1 - \delta_{\sigma,-1} \delta_{i,j}
\delta_{s,1}
         \right) \mbox{\textsf d}_s \A^+_{i1} \A^+_{js} \nonumber \\
  && \mbox{} - hh' \delta_{j,n} \delta_{t,m} \left[1 - \delta_{\sigma,-1} \left(
         \delta_{i,1} \delta_{s,1} + \delta_{i,1} \delta_{s,m} + \delta_{i,n}
         \delta_{s,1}\right)\right] d_i \mbox{\textsf d}_s \A^+_{11}
         \A^+_{is}\bigr\},  \label{eq:com1-h}
\end{eqnarray}
\begin{eqnarray}
  \left[\A_{is}, \A_{jt}\right\} & = & - (1 - \sigma P_{ij} P_{st})
         \bigl\{h \delta_{j,1} \left(1 - \delta_{\sigma,-1} \delta_{i,n}
\delta_{s,t}
         \right) d_i \A_{ns} \A_{it} \nonumber \\
  && \mbox{} + h' \delta_{t,1} \left(1 - \delta_{\sigma,-1} \delta_{i,j}
\delta_{s,m}
         \right) \mbox{\textsf d}_s \A_{im} \A_{js} \nonumber \\
  && \mbox{} + hh' \delta_{j,1} \delta_{t,1} \left[1 - \delta_{\sigma,-1} \left(
         \delta_{i,1} \delta_{s,m} + \delta_{i,n} \delta_{s,1} + \delta_{i,n}
         \delta_{s,m}\right)\right] d_i \mbox{\textsf d}_s \A_{nm}
\A_{is}\bigr\},
\end{eqnarray}
\begin{eqnarray}
  \left[\A_{is}, \A^+_{jt}\right\} & = & \delta_{i,j} \delta_{s,t}
         \left(\I + \sigma hh' d_i \mbox{\textsf d}_s \A^+_{11}
         \A_{nm}\right) \nonumber \\
  && \mbox{} + \sigma h \delta_{i,j} d_i \left[\A^+_{1t} \A_{ns}
         + h' \delta_{s,1} \delta_{t,m} \left(- \B_{1n} + h' \A^+_{11}
         \A_{nm}\right)\right] \nonumber \\
  && \mbox{} + \sigma h' \delta_{s,t} \mbox{\textsf d}_s \left[\A^+_{j1}
         \A_{im} + h \delta_{i,1} \delta_{j,n} \left(- \BB_{1m} + h \A^+_{11}
         \A_{nm}\right)\right] \nonumber \\
  && \mbox{} + \sigma h \delta_{i,1} \delta_{j,n} \left(- \BB_{ts}
         + h \A^+_{1t} \A_{ns}\right) + \sigma h' \delta_{s,1} \delta_{t,m}
         \left(- \B_{ji} + h' \A^+_{j1} \A_{im}\right) \nonumber \\
  && \mbox{} + \sigma hh' \delta_{i,1} \delta_{j,n} \delta_{s,1} \delta_{t,m}
         \left(\D - h \B_{1n} - h' \BB_{1m} + hh' \A^+_{11} \A_{nm}\right),
         \label{eq:com3-h}
\end{eqnarray}
where
\begin{equation}
  d_i = 2 - \delta_{i,1} - \delta_{i,n}, \qquad \mbox{\textsf d}_s = 2 -
\delta_{s,1}
  - \delta_{s,m},
\end{equation}
\begin{equation}
  \B_{ij} = \sum_u \mbox{\textsf d}_u \A^+_{iu} \A_{ju}, \qquad
  \BB_{st} = \sum_k d_k \A^+_{ks} \A_{kt}, \qquad
  \D = \sum_{ku} d_k \mbox{\textsf d}_u \A^+_{ku} \A_{ku},
\end{equation}
and $P_{ij}$ (resp.\ $P_{st}$) is the permutation operator acting on $i$,
$j$ (resp.\
$s$, $t$) indices.\par
%
%
In the $m=1$ case, Eqs.\ (\ref{eq:com1-h})--(\ref{eq:com3-h}) assume a much
simpler form
\begin{equation}
  \left[A^+_i, A^+_j\right\} = \left(1 - \sigma P_{ij}\right) \left[h
\delta_{j,n}
  \left(1 - \delta_{\sigma,-1} \delta_{i,1}\right) d_i A^+_1 A^+_i\right],
  \label{eq:com1-h1}
\end{equation}
\begin{equation}
  \left[A_i, A_j\right\} = - \left(1 - \sigma P_{ij}\right) \left[h \delta_{j,1}
  \left(1 - \delta_{\sigma,-1} \delta_{i,n}\right) d_i A_n A_i\right],
\end{equation}
\begin{equation}
  \left[A_i, A^+_j\right\} = \delta_{i,j} \left(I + \sigma h d_i A^+_1
A_n\right)
  + \sigma h \delta_{i,1} \delta_{j,n} \left(- \sum_k d_k A^+_k A_k + h
A^+_1 A_n
  \right).  \label{eq:com3-h1}
\end{equation}
\par
%
%
Let us next consider Eqs.\ (\ref{eq:compact1}), (\ref{eq:compact2-tilde}), and
(\ref{eq:compact3-tilde}), defining ${\cal A}^{(1)}_{q\sigma}(n,m)$ in the
$\left\{\A^{\prime +}_{is}, \tA'_{is}\right\}$ basis. Introducing transformed
$q$-bosonic (or $q$-fermionic) creation operators $\A^{\prime\prime +} =
\A^{\prime +} \g$ as before, and accordingly $\tA'' = \tA' \g$, we notice that
compatibility of the $\tA''$ and $\A''$ definitions with $\tA'' = \A''
\C''$, where $\C''
= C''_q {\cal C}''_{q^{\sigma}}$, leads to $C''_q = g^t C'_q g$, and ${\cal
C}''_{q^{\sigma}} = \mbox{\textsf g}^t {\cal C}'_{q^{\sigma}} \mbox{\textsf
g}$. A
simple calculation shows that for $n>1$
\begin{equation}
  C''_q = \sum_i (-1)^{n-i} q^{-(n-2i+1)/2} e_{ii'} + \eta
\left(q^{(n-1)/2} + (-1)^{n-1}
  q^{-(n-1)/2}\right) e_{nn},
\end{equation}
which can be rewritten as
\begin{eqnarray}
  C''_q & = & \sum_i (-1)^i q^{-(n-2i+1)/2} e_{ii'} + h \left(q^{(n-3)/2} +
q^{(n-5)/2}
          + \cdots + q^{-(n-1)/2}\right) e_{nn}, \nonumber \\
  & & \mbox{if\ } n = 2, 4, \ldots, \nonumber \\
  & = & \sum_i (-1)^{i-1} q^{-(n-2i+1)/2} e_{ii'} + \eta \left(q^{(n-1)/2}
          + q^{-(n-1)/2}\right) e_{nn}, \nonumber \\
  & & \mbox{if\ } n = 3, 5, \ldots.
\end{eqnarray}
We conclude that except for the trivial $n=1$ case, wherein we may set $C'_q =
C''_q = C_h = 1$, a contraction limit of $C''_q$ only exists for even $n$
values, and
is given by
\begin{equation}
  C_h \equiv \lim_{q\to1} C''_q = \sum_i (-1)^i e_{ii'} + (n-1) h e_{nn}.
\end{equation}
Similarly, for even $m$ values,
\begin{equation}
  {\cal C}_{h'} \equiv \lim_{q\to1} {\cal C}''_{q^{\sigma}} = \sum_s (-1)^s
e_{ss'}
  + (m-1) h' e_{mm}.
\end{equation}
\par
%
%
Restricting the range of $n$, $m$ values to $\{1, 2, 4, 6,\ldots\}$, we
obtain that
after transformation, Eqs.\ (\ref{eq:compact1}), (\ref{eq:compact2-tilde}), and
(\ref{eq:compact3-tilde}) contract into
\begin{equation}
  \A^+_1 \A^+_2 = \sigma \A^+_2 \A^+_1 R_h {\cal R}_{h'},
  \label{eq:compact1-hbis}
\end{equation}
\begin{equation}
  \tA_1 \tA_2 = \sigma \tA_2 \tA_1 R_h {\cal R}_{h'},
\label{eq:compact2-tilde-h}
\end{equation}
\begin{equation}
  \tA_2 \A^+_1 = \C_{12} + \sigma \A^+_1 \tA_2 \tilde{R}^{-1}_h \tilde{{\cal
  R}}^{-1}_{h'}, \label{eq:compact3-tilde-h}
\end{equation}
where $\C = C_h {\cal C}_{h'}$,
\begin{eqnarray}
  \tilde{R}_h & \equiv & \lim_{q\to1} \left(g^{-1} \otimes g^{-1}\right)
\tilde{R}'_q
          (g \otimes g) \nonumber \\
  & = & C^{-1}_{h,1} \left(R^{-1}_h\right)^{t_1} C_{h,1} = C^{-1}_{h,2}
          \left(R^{t_2}_h\right)^{-1} C_{h,2} \nonumber \\
  & = & \sum_{ij} e_{ii} \otimes e_{jj} - h \sum_i (-1)^i d_i \left(e_{1i}
\otimes
          e_{1i'} + e_{in} \otimes e_{i'n}\right) + (2n-3) h^2 e_{1n}
\otimes e_{1n},
\end{eqnarray}
and $\tilde{{\cal R}}_{h'}$ is defined in the same way.\par
%
%
Again the same procedure applied to Eqs.\ (\ref{eq:compact1}),
(\ref{eq:compact2-tilde}), and (\ref{eq:compact3bis-tilde}), defining ${\cal
A}^{(2)}_{q\sigma}(n,m)$ in the $\left\{\A^{\prime +}_{is},
\tA'_{is}\right\}$ basis, leads to Eqs.\
(\ref{eq:compact1-hbis})--(\ref{eq:compact3-tilde-h}), already obtained for
${\cal
A}^{(1)}_{q\sigma}(n,m)$. We conclude that for $n,m \in \{1, 2, 4,
6,\ldots\}$, such equations yield another form of the $\mbox{GL}_h(n) \times
\mbox{GL}_{h'}(m)$-covariant ($hh'$)-bosonic (or ($hh'$)-fermionic) algebra
${\cal
A}_{hh'\sigma}(n,m)$, defined in Eqs.\
(\ref{eq:compact1-h})--(\ref{eq:compact3-h})
for arbitrary $n$, $m$ values. The counterpart of transformation
(\ref{eq:phi}) is
now
\begin{equation}
  \varphi\left(\A^+\right) = \A^+ T {\cal T}, \qquad
  \varphi\left(\tA\right) = \tA \tilde{T} \tilde{{\cal T}},
\end{equation}
where $T_{ij} \in \mbox{GL}_h(n)$, ${\cal T}_{st} \in \mbox{GL}_{h'}(m)$,
$\tilde{T} = C^{-1}_h \left(T^{-1}\right)^t C_h$, and $\tilde{{\cal T}} = {\cal
C}^{-1}_{h'} \left({\cal T}^{-1}\right)^t {\cal C}_{h'}$. However, for $n$
and/or $m
\in \{3, 5, 7,\ldots\}$, the contraction procedure does not preserve the
equivalence between the two forms of ${\cal A}^{(1)}_{q\sigma}(n,m)$ or
${\cal A}^{(2)}_{q\sigma}(n,m)$, corresponding to the $\left\{\A^{\prime
+}_{is}, \A'_{is}\right\}$ and $\left\{\A^{\prime +}_{is},
\tA'_{is}\right\}$ bases,
respectively, since only the former has a limit. It should be stressed that such
results are entirely new, since Fiore (1997) did not consider any $\tA'_{is}$
operators.\par
%
%
In componentwise form, Eq.\ (\ref{eq:compact1-hbis}) becomes Eq.\
(\ref{eq:com1-h}), Eq.\ (\ref{eq:compact2-tilde-h}) assumes a similar form,
while
Eq.\ (\ref{eq:compact3-tilde-h}) leads to the following relation
\begin{eqnarray}
  \left[\tA_{is}, \A^+_{jt}\right\} & = & \delta_{i',j} \delta_{s',t}
         (-1)^{i+s} \left(\I + \sigma hh' d_i \mbox{\textsf d}_s \A^+_{11}
         \tA_{11}\right) \nonumber \\
  & & \mbox{} - \delta_{i',j} (-1)^i \bigl\{\sigma h d_i \A^+_{1t}
        \tA_{1s} + h' \delta_{s,m} \delta_{t,m} \bigl[(m-1) \I
        \nonumber \\
  & & \mbox{} + \sigma h d_i \tB_{11} + \sigma (2m-3) hh' d_i
        \A^+_{11} \tA_{11}\bigr]\bigr\} \nonumber \\
  & & \mbox{} - \delta_{s',t} (-1)^s \bigl\{\sigma h' \mbox{\textsf d}_s
        \A^+_{j1} \tA_{i1} + h \delta_{i,n} \delta_{j,n}
        \bigl[(n-1) \I \nonumber \\
  & & \mbox{} + \sigma h' \mbox{\textsf d}_s \tBB_{11}
        + \sigma (2n-3) hh' {\textsf d}_s \A^+_{11}
        \tA_{11}\bigr]\bigr\} \nonumber \\
  & & \mbox{} + \sigma h \delta_{i,n} \delta_{j,n} \left[\tBB_{ts}
        + (2n-3) h \A^+_{1t} \tA_{1s}\right] \nonumber \\
  & & \mbox{} + \sigma h' \delta_{s,m} \delta_{t,m} \left[\tB_{ji}
        + (2m-3) h' \A^+_{j1} \tA_{i1}\right] \nonumber \\
  & & \mbox{} + hh' \delta_{i,n} \delta_{j,n} \delta_{s,m} \delta_{t,m}
        \bigl[(n-1)(m-1) \I + \sigma \tD + \sigma (2n-3) h
        \tB_{11} \nonumber \\
  & & {} + \sigma (2m-3) h' \tBB_{11}
        + \sigma (2n-3) (2m-3) hh' \A^+_{11} \tA_{11}\bigr],
        \label{eq:com3tilde-h}
\end{eqnarray}
where
\begin{equation}
  \tB_{ij} = \sum_u (-1)^u \mbox{\textsf d}_u \A^+_{iu}
  \tA_{ju'}, \quad  \tBB_{st} = \sum_k (-1)^k d_k
  \A^+_{ks} \tA_{k't}, \quad \tD = \sum_{ku}
  (-1)^{k+u} d_k \mbox{\textsf d}_u \A^+_{ku} \tA_{k'u'}.
\end{equation}
\par
%
%
In the $m=1$ case, Eq.\ (\ref{eq:com3tilde-h}) assumes the simpler form
\begin{eqnarray}
  \left[\tilde{A}_i, A^+_j\right\} & = & \delta_{i',j} (-1)^{i+1} \left(I +
\sigma h d_i
           A^+_1 \tilde{A}_1\right) \nonumber \\
  && \mbox{} + h \delta_{i,n} \delta_{j,n} \Bigl[(n-1) I + \sigma \sum_k
(-1)^k d_k
           A^+_k \tilde{A}_{k'} + \sigma (2n-3) h A^+_1 \tilde{A}_1\Bigr],
           \label{eq:com3tilde-h1}
\end{eqnarray}
where $\tilde{A} = A C_h$.\par
%
%
In the next section, by making explicit use of the U$_h$(sl(2)) CGC
determined by
Van der Jeugt (1998), we plan to show that whenever $n=2$, and $m=1$ or 2, the
(anti)commutators (\ref{eq:com1-h}) and (\ref{eq:com3tilde-h}) can be rewritten
in coupled form as in the $q$-deformed case.\par
%
%
\section{\boldmath SPECIAL CASES $n=2$, $m=1$ AND $n=m=2$}
\label{sec:special}
\setcounter{equation}{0}
Let us first consider the $n=2$, $m=1$ case, wherein
\begin{equation}
  R_h = \left(\begin{array}{cccc}
        1 & h & -h & h^2 \\[0.1cm]
        0 & 1 & 0  & h \\[0.1cm]
        0 & 0 & 1  & -h \\[0.1cm]
        0 & 0 & 0  & 1
        \end{array}\right), \qquad
  C_h = \left(\begin{array}{cc}
        0 & -1 \\[0.1cm]
        1 & h
        \end{array}\right), \qquad
  {\cal R}_{h'} = {\cal C}_{h'} = 1.  \label{eq:special}
\end{equation}
\par
%
%
{}From Eqs.\ (\ref{eq:com1-h1})--(\ref{eq:com3-h1}), and
(\ref{eq:com3tilde-h1}),
it follows that the defining relations of the GL$_h$(2)-covariant $h$-bosonic
algebra ${\cal A}_{h+}(2,1)$ are given by
\begin{equation}
  \left[A^+_1, A^+_2\right] = h \left(A^+_1\right)^2, \qquad \left[A_1,
A_2\right]
  = h A^2_2,  \label{eq:com-h1-b}
\end{equation}
\begin{equation}
  \left[A_2, A^+_1\right] = 0, \qquad \left[A_1, A^+_2\right] = h \left(-
A^+_1 A_1
  - A^+_2 A_2 + h A^+_1 A_2\right),
\end{equation}
\begin{equation}
  \left[A_1, A^+_1\right] = \left[A_2, A^+_2\right] = I + h A^+_1 A_2,
\end{equation}
in the $\left\{A^+_1, A^+_2, A_1, A_2\right\}$ basis, and by
\begin{equation}
  \left[A^+_1, A^+_2\right] = h \left(A^+_1\right)^2, \qquad \bigl[\tilde{A}_1,
        \tilde{A}_2\bigr] = h \tilde{A}^2_1, \label{eq:comtilde-h1-b1}
\end{equation}
\begin{equation}
  \bigl[\tilde{A}_1, A^+_1\bigr] = 0, \qquad \bigl[\tilde{A}_2, A^+_2\bigr] =
  h \bigl(I - A^+_1 \tilde{A}_2 + A^+_2 \tilde{A}_1 + h A^+_1 \tilde{A}_1\bigr),
\end{equation}
\begin{equation}
  \bigl[\tilde{A}_1, A^+_2\bigr] = - \bigl[\tilde{A}_2, A^+_1\bigr] = I + h
A^+_1
  \tilde{A}_1, \label{eq:comtilde-h1-b3}
\end{equation}
in the $\left\{A^+_1, A^+_2, \tilde{A}_1, \tilde{A}_2\right\}$ one.\par
%
%
Similarly, for the $h$-fermionic algebra ${\cal A}_{h-}(2,1)$, we obtain
\begin{equation}
  \left\{A^+_1, A^+_1\right\} = \left\{A^+_1, A^+_2\right\} = 0, \qquad
  \left\{A^+_2, A^+_2\right\} = 2h A^+_1 A^+_2,
\end{equation}
\begin{equation}
  \left\{A_1, A_1\right\} = 2h A_1 A_2, \qquad \left\{A_1, A_2\right\} =
  \left\{A_2, A_2\right\} = 0,
\end{equation}
\begin{equation}
  \left\{A_2, A^+_1\right\} = 0, \qquad \left\{A_1, A^+_2\right\} = h
\left(A^+_1
  A_1 + A^+_2 A_2 - h A^+_1 A_2\right),
\end{equation}
\begin{equation}
  \left\{A_1, A^+_1\right\} = \left\{A_2, A^+_2\right\} = I - h A^+_1 A_2,
\end{equation}
and
\begin{equation}
  \left\{A^+_1, A^+_1\right\} = \left\{A^+_1, A^+_2\right\} = 0, \qquad
  \left\{A^+_2, A^+_2\right\} = 2h A^+_1 A^+_2,  \label{eq:comtilde-h1-f1}
\end{equation}
\begin{equation}
  \left\{\tilde{A}_1, \tilde{A}_1\right\} = \left\{\tilde{A}_1,
\tilde{A}_2\right\} =
  0, \qquad \left\{\tilde{A}_2, \tilde{A}_2\right\} = 2h \tilde{A}_1
\tilde{A}_2,
\end{equation}
\begin{equation}
  \bigl\{\tilde{A}_1, A^+_1\bigr\} = 0, \qquad \bigl\{\tilde{A}_2,
A^+_2\bigr\} =
  h \bigl(I + A^+_1 \tilde{A}_2 - A^+_2 \tilde{A}_1 - h A^+_1 \tilde{A}_1\bigr),
\end{equation}
\begin{equation}
  \bigl\{\tilde{A}_1, A^+_2\bigr\} = - \bigl\{\tilde{A}_2, A^+_1\bigr\} = I
- h A^+_1
  \tilde{A}_1,  \label{eq:comtilde-h1-f4}
\end{equation}
respectively.\par
%
%
The operators $\left(A^+_1, A^+_2\right)$, and $\left(\tilde{A}_1,
\tilde{A}_2\right)$ may be considered as the components $m=1/2$, and $m=-1/2$
of ITO of rank 1/2, or spinors, with respect to the quantum algebra
U$_h$(sl(2)). By
using a nonlinear invertible map between the generators of U$_h$(sl(2)) and
U(sl(2)) (Abdesselam {\em et.al.}, 1996), and considering the adjoint
action of the
former on such spinors, Aizawa (1998) recently realized them in terms of
standard bosonic or fermionic operators $a^+_1$,
$a^+_2$, $a_1$, $a_2$. For the standard form of sl(2) generators
\begin{equation}
   J_+ = a^+_1 a_2, \qquad J_- = a^+_2 a_1, \qquad J_0 = \case{1}{2}
\left(a^+_1 a_1
  - a^+_2 a_2\right),
\end{equation}
the realizations read$^1$
\begin{equation}
  A^+_1 = \left(1 - \case{h}{2} J_+\right)^{-1} a^+_1, \qquad A^+_2 = \left(1 -
  \case{h}{2} J_+\right) a^+_2 + \case{h}{2} \left(A^+_1 - 2 a^+_1 J_0\right),
  \label{eq:realization-b1}
\end{equation}
\begin{equation}
  \tilde{A}_1= \left(1 - \case{h}{2} J_+\right)^{-1} a_2, \qquad \tilde{A}_2
  = - \left(1 - \case{h}{2} J_+\right) a_1 + \case{h}{2} \left(\tilde{A}_1
- 2 a_2
  J_0\right),
  \label{eq:realization-b2}
\end{equation}
in the $h$-bosonic case, and
\begin{equation}
  A^+_1 = a^+_1, \qquad A^+_2 = a^+_2 - 2h a^+_1 J_0,  \label{eq:realization-f1}
\end{equation}
\begin{equation}
  \tilde{A}_1= a_2, \qquad \tilde{A}_2 = - a_1 - 2h a_2 J_0,
  \label{eq:realization-f2}
\end{equation}
in the $h$-fermionic one. As expected, the operators (\ref{eq:realization-b1}),
(\ref{eq:realization-b2}) and (\ref{eq:realization-f1}),
(\ref{eq:realization-f2})
satisfy Eqs.\ (\ref{eq:comtilde-h1-b1})--(\ref{eq:comtilde-h1-b3}) and
(\ref{eq:comtilde-h1-f1})--(\ref{eq:comtilde-h1-f4}), respectively.\par
%
%
Let us now introduce coupled (anti)commutators, defined as in Eq.\
(\ref{eq:coupledcom}) by
\begin{equation}
  \left[T^{j_1}, U^{j_2}\right\}^J_M = \left[T^{j_1} \times U^{j_2}\right]^J_M
  - \sigma (-1)^{\epsilon} \left[U^{j_2} \times T^{j_1}\right]^J_M.
\end{equation}
Here $T^{j_1}$ and $U^{j_2}$ denote two ITO of rank $j_1$ and $j_2$ with
respect to U$_h$(sl(2)), respectively, $\epsilon$ is defined as in
Eqs.~(\ref{eq:epsilon}), (\ref{eq:phi-lambda}) by $\epsilon = j_1 + j_2 -
J$, and
\begin{equation}
  \left[T^{j_1} \times U^{j_2}\right]^J_M = \sum_{m_1m_2} \langle j_1
  m_1, j_2 m_2 | J M \rangle_h\, T^{j_1}_{m_1} U^{j_2}_{m_2},
\end{equation}
where $\langle\, ,\, | \,\rangle_h$ denotes a U$_h$(sl(2)) CGC (Van der
Jeugt, 1998).
The values of the latter needed for coupling spinors are given in Table~I.
By using
them, Eqs.\ (\ref{eq:comtilde-h1-b1})--(\ref{eq:comtilde-h1-b3}) and
(\ref{eq:comtilde-h1-f1})--(\ref{eq:comtilde-h1-f4}) can be recast in the
compact
forms
\begin{equation}
  \bigl[A^+, A^+\bigr]^0_0 = \bigl[\tilde{A}, \tilde{A}\bigr]^0_0 =
\bigl[\tilde{A},
  A^+\bigr]^1_M = 0, \qquad \bigl[\tilde{A}, A^+\bigr]^0_0 = \sqrt{2}\, I,
  \label{eq:coupledcom-h-b1}
\end{equation}
and
\begin{equation}
  \bigl\{A^+, A^+\bigr\}^1_M = \bigl\{\tilde{A}, \tilde{A}\bigr\}^1_M
  = \bigl\{\tilde{A}, A^+\bigr\}^1_M = 0, \qquad \bigl\{\tilde{A},
A^+\bigr\}^0_0
  = \sqrt{2}\, I, \label{eq:coupledcom-h-f1}
\end{equation}
respectively.\par
%
%
Let us next consider the $n=m=2$ case, wherein ${\cal R}_{h'}$ and ${\cal
C}_{h'}$
are defined as $R_h$ and $C_h$ in Eq.\ (\ref{eq:special}). Relations
similar to Eqs.\
(\ref{eq:com-h1-b})--(\ref{eq:comtilde-h1-f4}) can be easily written. The
operators
$\A^+_{is}$ ($i$, $s=1$, 2), and $\tA_{is}$ ($i$, $s=1$, 2) may now be
considered as
the components of double spinors with respect to $\mbox{U}_h(\mbox{sl}(2))
\times
\mbox{U}_{h'}(\mbox{sl}(2))$. Defining coupled (anti)commutators by
\begin{equation}
  \left[T^{j_1j'_1}, U^{j_2j'_2}\right\}^{JJ'}_{MM'} = \left[T^{j_1j'_1} \times
  U^{j_2j'_2}\right]^{JJ'}_{MM'} - \sigma (-1)^{\epsilon} \left[U^{j_2j'_2}
\times
  T^{j_1j'_1}\right]^{JJ'}_{MM'},
\end{equation}
where $\epsilon = j_1 + j_2 - J + j'_1 + j'_2 - J'$, and
\begin{equation}
  \left[T^{j_1j'_1} \times U^{j_2j'_2}\right]^{JJ'}_{MM'} =
\sum_{m_1m_2m'_1m'_2}
  \langle j_1 m_1, j_2 m_2 | J M \rangle_h\, \langle j'_1 m'_1, j'_2 m'_2 |
J' M'
  \rangle_{h'}\, T^{j_1j'_1}_{m_1m'_1} U^{j_2j'_2}_{m_2m'_2},
\end{equation}
we easily obtain that the double spinors $\A^+$ and $\tA$ satisfy the relations
\begin{equation}
  \bigl[\A^+, \A^+\bigr]^{10}_{M0} = \bigl[\A^+, \A^+\bigr]^{01}_{0M'} =
  \bigl[\tA, \tA\bigr]^{10}_{M0} = \bigl[\tA, \tA\bigr]^{01}_{0M'} = 0,
  \label{eq:coupledcom-h-b2-1}
\end{equation}
\begin{equation}
  \bigl[\tA, \A^+\bigr]^{JJ'}_{MM'}= 2 \delta_{J,0} \delta_{J',0} \delta_{M,0}
  \delta_{M',0} \I,  \label{eq:coupledcom-h-b2-2}
\end{equation}
and
\begin{equation}
  \bigl\{\A^+, \A^+\bigr\}^{11}_{MM'} = \bigl\{\A^+, \A^+\bigr\}^{00}_{00} =
  \bigl\{\tA, \tA\bigr\}^{11}_{MM'} = \bigl\{\tA, \tA\bigr\}^{00}_{00} = 0,
  \label{eq:coupledcom-h-f2-1}
\end{equation}
\begin{equation}
  \bigl\{\tA, \A^+\bigr\}^{JJ'}_{MM'}= 2 \delta_{J,0} \delta_{J',0} \delta_{M,0}
  \delta_{M',0} \I,  \label{eq:coupledcom-h-f2-2}
\end{equation}
in the ($hh'$)-bosonic and ($hh'$)-fermionic cases, respectively.\par
%
%
It is remarkable that Eqs.\ (\ref{eq:coupledcom-h-b1}) (resp.\
(\ref{eq:coupledcom-h-f1}), and (\ref{eq:coupledcom-h-b2-1}),
(\ref{eq:coupledcom-h-b2-2}) (resp.\ (\ref{eq:coupledcom-h-f2-1}),
(\ref{eq:coupledcom-h-f2-2})) are formally identical with those for bosonic
(resp.\
fermionic) ITO with respect to the Lie algebras sl(2) and
$\mbox{sl}(2) \times \mbox{sl}(2)$, respectively. Contrary to what happens
in the
$q$-bosonic (or $q$-fermionic) case where the (anti)commutators are
$q$-deformed (see Eqs.~(\ref{eq:coupledcom3-1}) and (\ref{eq:coupledcom3-2})),
here all the dependence upon the deforming parameters $h$, $h'$ is
contained in the
coupling coefficients.\par
%
%
\section{CONCLUSION}
\label{sec:conclusion}
In the present paper, we did show that the contraction technique,
previously used to
construct Jordanian deformations of Lie groups from standard
ones (Aghamohammadi {\em et al.}, 1995; Alishahiha, 1995) can be applied to the
$\mbox{GL}_q(n) \times \mbox{GL}_q(m)$-covariant $q$-bosonic (or
$\mbox{GL}_q(n) \times \mbox{GL}_{q^{-1}}(m)$-covariant $q$-fermionic) algebras
${\cal A}^{(\alpha)}_{q\pm}(n,m)$, $\alpha = 1, 2$ (Quesne, 1993; Quesne,
1994; Fiore, 1998), to yield some $\mbox{GL}_h(n) \times
\mbox{GL}_{h'}(m)$-covariant ($hh'$)-bosonic (or ($hh'$)-fermionic)
algebras ${\cal
A}_{hh'\pm}(n,m)$. In this process, the arbitrariness present in the
$q$-deformed
case disappears as the algebras ${\cal A}^{(1)}_{q\pm}(n,m)$ and ${\cal
A}^{(2)}_{q\pm}(n,m)$ have the same contraction limit ${\cal
A}_{hh'\pm}(n,m)$.\par
%
%
When using a basis $\{\A^{\prime+}_{is}, \A'_{is}\}$ of ${\cal
A}^{(\alpha)}_{q\pm}(n,m)$, wherein the annihilation operators $\A'_{is}$ are
contragredient to the creation ones $\A^{\prime+}_{is}$, this contraction
procedure
can be carried out for any $n$, $m$ values. The resulting defining relations of
${\cal A}_{hh'\pm}(n,m)$ were written in the contracted basis $\{\A^+_{is},
\A_{is}\}$, both in compact form in terms of the defining $R_h$ and ${\cal
R}_{h'}$-matrices of GL$_h$($n$) and GL$_{h'}$($m$), respectively, and in
componentwise form. They may be considered as a special case of the defining
relations of $\cal H$-covariant deformed bosonic (or fermionic) algebras for
triangular Hopf algebras $\cal H$, recently obtained by Fiore (1997) by another
procedure.\par
%
%
When using instead a basis $\{\A^{\prime+}_{is}, \tA'_{is}\}$ of ${\cal
A}^{(\alpha)}_{q\pm}(n,m)$, wherein the annihilation operators $\tA'_{is}$
are ITO
with respect to the quantum algebra $\mbox{U}_q(\mbox{gl}(n)) \times
 \mbox{U}_{q^{\pm1}}(\mbox{gl}(m))$, we obtained some new and interesting
results. We did indeed establish that in such a case a contraction limit
only exists
whenever $n$, $m \in \{1, 2, 4, 6, \ldots\}$, hence showing that for $n$
and/or $m
\in \{3, 5, 7, \ldots\}$, the contraction procedure does not preserve the
equivalence between the two forms of ${\cal A}^{(\alpha)}_{q\pm}(n,m)$,
corresponding to the $\{\A^{\prime+}_{is}, \A'_{is}\}$ and $\{\A^{\prime+}_{is},
\tA'_{is}\}$ bases. When a limit does exist, the defining relations of
${\cal A}_{hh'\pm}(n,m)$ were written in the contracted basis $\{\A^+_{is},
\tA_{is}\}$, both in compact form in terms of $R_h$ and ${\cal R}_{h'}$, and in
componentwise form.\par
%
%
Such a basis is essential to express the defining relations of ${\cal
A}_{hh'\pm}(n,m)$ in another compact form in terms of coupled (anti)commutators,
thereby enhancing the transformation properties of the generators under the
quantum algebra dual to $\mbox{GL}_h(n) \times \mbox{GL}_{h'}(m)$. We did prove
this point in the $n=2$, $m=1$, and $n=m=2$ cases, where the dual quantum
algebras are known, and the U$_h$(sl(2)) CGC determined by Van der Jeugt (1998)
can be used. Furthermore, we did check that the $h$-bosonic and
$h$-fermionic ITO
of rank 1/2 with respect to U$_h$(sl(2)), constructed by Aizawa (1998), satisfy
the defining relations of ${\cal A}_{h\pm}(2,1)$. From the examples
considered, we
concluded that the algebras ${\cal A}_{hh'\pm}(n,m)$ are much closer to the
standard Heisenberg (or Clifford) algebras ${\cal A}_{\pm}(n,m)$ than the
$q$-deformed ones, ${\cal A}^{(\alpha)}_{q\pm}(n,m)$. This may be an
advantage in
some physical applications.\par
%
%
\section*{ACKNOWLEDGMENT}
The author is a Research Director of the National Fund for Scientific Research
(FNRS), Belgium.\par
%
%
\newpage
\setlength{\parindent}{0cm}
\section*{FOOTNOTES}

$^1$The realization of sl(2) used in Eqs.\ (\ref{eq:realization-f1}) and
(\ref{eq:realization-f2}) differs from that considered by Aizawa (1998).
There are
also some changes of phase conventions with respect to the same reference in
Eqs.\ (\ref{eq:realization-b1})--(\ref{eq:realization-f2}).\par
%
%
\newpage
\section*{REFERENCES}

Abdesselam, B., Chakrabarti, A., and Chakrabarti, R.\ (1996).\ {\em Modern
Physics
Letters} A, {\bf 11}, 2883.

Abdesselam, B., Chakrabarti, A., and Chakrabarti, R.\ (1997).\ {\em
International
Journal of Modern Physics} A, {\bf 12}, 2301.

Abdesselam, B., Chakrabarti, A., and Chakrabarti, R.\ (1998a).\ {\em Modern
Physics
Letters} A, {\bf 13}, 779.

Abdesselam, B., Chakrabarti, A., Chakrabarti, R., and Segar, J.\ (1998b).\
Maps and
twists relating U(sl(2)) and the nonstandard U$_h$(sl(2)): unified construction.
Preprint math.QA/9807100.

Aghamohammadi, A.\ (1993).\ {\em Modern Physics Letters} A, {\bf 8}, 2607.

Aghamohammadi, A., Khorrami, M., and Shariati, A.\ (1995).\ {\em Journal of
Physics} A, {\bf 28}, L225.

Aizawa, N.\ (1997).\ {\em Journal of Physics} A, {\bf 30}, 5981.

Aizawa, N.\ (1998).\ {\em Journal of Physics} A, {\bf 31}, 5467.

Alishahiha, M.\ (1995).\ {\em Journal of Physics} A, {\bf 28}, 6187.

Aneva, B. L., Dobrev, V. K., and Milov, S. G.\ (1997).\ {\em Journal of
Physics} A,
{\bf 30}, 6769.

Ballesteros, A., and Herranz, F. J.\ (1996).\ {\em Journal of Physics} A,
{\bf 29},
L311.

Biedenharn, L. C.\ (1989).\ {\em Journal of Physics} A, {\bf 22}, L873.

Biedenharn, L. C.\ (1990).\ A $q$-boson realization of the quantum group
SU$_q$(2)
and the theory of $q$-tensor operators, in {\em Quantum Groups, Proceedings
of the
8th Workshop on Mathematical Physics (Clausthal, FRG, 1989)}, H.\ D.\
Doebner and
J.\ D.\ Hennig, eds., Lecture Notes in Physics, Vol.\ {\bf 370}, Springer,
Berlin.

Chaichian, M., and Kulish, P.\ (1990).\ {\em Physics Letters} B, {\bf 234}, 72.

Chakrabarti, R., and Quesne, C.\ (1998).\ On Jordanian U$_{h,\alpha}$(gl(2))
algebra and its $T$ matrices via a contraction method. Universit\'e Libre de
Bruxelles preprint ULB/229/CQ/98/3, math.QA/9811064.

Demidov, E. E., Manin, Yu. I., Mukhin, E. E., and Zhdanovich, D. V.\
(1990).\ {\em
Progress in Theoretical Physics Supplement}, {\bf 102}, 203.

Dobrev, V. K.\ (1996).\ In {\em Proceedings of the 10th International Conference
`Problems of Quantum Field Theory', (Alushta, Crimea, Ukraine, 13--18.5.1996)},
D.\ Shirkov, D.\ Kazakov, and A.\ Vladimirov, eds., JINR E2-96-369, Dubna.

Drinfeld, V. G.\ (1987).\ Quantum groups, in {\em Proceedings of the
International
Congress of Mathematicians, Berkeley, California, 1986},  A.\ M.\ Gleason, ed.,
American Mathematical Society, Providence.

Fiore, G.\ (1997).\ On Bose-Fermi statistics, quantum group symmetry, and second
quantization, in {\em Proceedings of the Quantum Group Symposium at Group21
(Goslar, Germany, 1996)}, H.-D.\ Doebner and V.\ K.\ Dobrev, eds., Heron
Press, Sofia.

Fiore, G.\ (1998).\ {\em Journal of Physics} A {\bf 31}, 5289.

Hayashi, T.\ (1990).\ {\em Communications in Mathematical Physics}, {\bf 127},
129.

Kupershmidt, B. A.\ (1992).\ {\em Journal of Physics} A, {\bf 25}, L1239.

Macfarlane, A. J.\ (1989).\ {\em Journal of Physics} A, {\bf 22}, 4581.

Majid, S.\ (1990).\ {\em International Journal of Modern Physics} A, {\bf 5}, 1.

Ohn, C.\ (1992).\ {\em Letters in Mathematical Physics}, {\bf 25}, 85.

Parashar, P.\ (1998).\ {\em Letters in Mathematical Physics}, {\bf 45}, 105.

Pusz, W.\ (1989).\ {\em Reports on Mathematical Physics}, {\bf 27}, 349.

Pusz, W., and Woronowicz, S. L.\ (1989).\ {\em Reports on Mathematical Physics},
{\bf27}, 231.

Quesne, C.\ (1993). {\em Physics Letters} B, {\bf 298}, 344.

Quesne, C.\ (1994). {\em Physics Letters} B, {\bf 322}, 344.

Shariati, A., Aghamohammadi, A., and Khorrami, M.\ (1996). {\em Modern Physics
Letters} A, {\bf 11}, 187.

Sun, C.-P., and Fu, H.-C.\ (1989).\ {\em Journal of Physics} A, {\bf 22}, L983.

Van der Jeugt, J.\ (1998).\ {\em Journal of Physics} A, {\bf 31}, 1495.

Zakrzewski, S.\ (1991).\ {\em Letters in Mathematical Physics}, {\bf 22}, 287.
%
%
\newpage

\begin{table}[h]

\caption{Values of U$_h$(sl(2)) CGC $\left\langle \case{1}{2} m_1, \case{1}{2}
m_2 | J M \right\rangle_h$.}

\vspace{1cm}
\begin{tabular}{lcccc}
  \hline\\[-0.2cm]
  & $J=M=1$ & $J=1$, $M=0$ & $J=-M=1$ & $J=M=0$\\[0.2cm]
  \hline\\[-0.2cm]
  $m_1=m_2=1/2$ & 1 & 0 & $(h/2)^2$ & $- h/\sqrt{2}$\\[0.2cm]
  $m_1=-m_2=1/2$ & 0 & $1/\sqrt{2}$ & $- h/2$ & $1/\sqrt{2}$\\[0.2cm]
  $m_1=-m_2=-1/2$ & 0 & $1/\sqrt{2}$ & $h/2$ & $- 1/\sqrt{2}$\\[0.2cm]
  $m_1=m_2=-1/2$ & 0 & 0 & 1 & 0\\[0.2cm]
  \hline
\end{tabular}

\end{table}

\end{document}